\documentclass{amsart}
\usepackage{latexsym}
\usepackage{amsmath, amssymb}
\usepackage[all]{xy}
\newtheorem{dfs}{Definition}[section]
\newtheorem{lms}[dfs]{Lemma}
\newtheorem{thms}[dfs]{Theorem}
\newtheorem{props}[dfs]{Proposition}

\begin{document}

\title[Non-isomorphic C$^*$-algebras with identical $\mathrm{K}$-theory]{An infinite family 
of non-isomorphic C$^*$-algebras with identical $\mathrm{K}$-theory}

\author{Andrew S. Toms}

\maketitle

\begin{abstract}
We exhibit a countably infinite family of simple, separable, nuclear, and mutually non-isomorphic C$^*$-algebras which agree
on $\mathrm{K}$-theory and traces.  The algebras do not absorb the Jiang-Su algebra $\mathcal{Z}$ tensorially,
answering a question of N. C. Phillips.  They are also pairwise shape and Morita equivalent,
confirming a conjecture from our earlier work.  
The distinguishing invariant is the radius of comparison, a non-stable invariant of the Cuntz semigroup.
\end{abstract}

\section{Introduction}\label{intro}

In 1997 J. Villadsen gave an example of a simple, separable, and nuclear C$^*$-algebra whose
ordered $\mathrm{K}_0$-group exhibited pathological behaviour, answering a long-standing open
question.  Since his discovery, those working on the structure of nuclear C$^*$-algebras have 
come to realise that there is, apparently, a dichotomy in their field:  algebras are
either amenable to classification via $\mathrm{K}$-theory, or so wild as to resist most attempts
to organise them.  Villadsen's example was simply the first of the ``wild'' algebras. 

There are several proposed characterisations of those simple, separable, and nuclear C$^*$-algebras
which will prove amenable to classification.  The property of absorbing the Jiang-Su algebra $\mathcal{Z}$
tensorially --- being {\it $\mathcal{Z}$-stable} --- is chief among them.  The main reason for this
is the existence of several examples due first to R{\o}rdam (\cite{R3}) and later the author 
(\cite{To1}, \cite{To2}) exhibiting pairs of simple, separable, and nuclear C$^*$-algebras which agree
on $\mathrm{K}$-theory and traces but are non-isomorphic.  The algebras are, in each case, distinguished 
by the fact that one of them is not $\mathcal{Z}$-stable.  They leave open the following question, 
posed by N. Christopher Phillips: 
\begin{quotation} \normalsize
Are there simple, separable, nuclear, and non-$\mathcal{Z}$-stable C$^*$-algebras 
which agree on $\mathrm{K}$-theory and traces, yet are non-isomorphic? 
\end{quotation}
\vspace{1mm}
The idea is this:  eliminate the possibility that the failure of
$\mathrm{K}$-theory to be a complete invariant for simple, separable, and nuclear C$^*$-algebras
is due simply to the need to differentiate between $\mathcal{Z}$-stable and non-$\mathcal{Z}$-stable
algebras.  Phillips' question is a specific instance of a more general question, namely, to what extent
can separable, nuclear, and non-isomorphic C$^*$-algebras be similar?  In the sequel 
we give a positive answer to Phillips' question,
and make some interesting progress on the more general question, too.

The major difficulty in answering Phillips' question lies in the fact that non-$\mathcal{Z}$-stable
C$^*$-algebras are difficult to control.  We cannot, to date, elicit prescribed $\mathrm{K}$-theoretic
data from them.  And if one manages to produce two such algebras with the same $\mathrm{K}$-theory
and tracial state space, how can one tell them apart?  In the sequel we overcome these challenges 
by studying a non-stable invariant of the Cuntz semigroup --- the {\it radius of comparison} --- introduced
in \cite{To3}.  We obtain:

\vspace{2mm}
\noindent
{\bf Theorem.} {\it
There is a simple, separable, nuclear, and 
non-$\mathcal{Z}$-stable C$^*$-algebra $A$ such that for any natural numbers $m \neq n$ we have:
\begin{enumerate}
\item[(i)] $\mathrm{M}_n(A) \ncong \mathrm{M}_m(A)$;
\item[(ii)] $\mathrm{M}_n(A)$ and $\mathrm{M}_m(A)$ are shape and (evidently) Morita equivalent;
\item[(iii)] $\mathrm{M}_n(A)$ and $\mathrm{M}_m(A)$ have the same scaled ordered $\mathrm{K}$-theory and tracial state space.
\end{enumerate} }
\noindent
The algebras of the theorem thus agree 
on stable isomorphism invariants, and on homotopy invariant 
functors which commute with sequential inductive limits.  They moreover agree on every version
of non-commutative dimension for C$^*$-algebras:  the real, stable, tracial topological, and
decomposition ranks.    The Morita equivalence of simple, separable,
nuclear, and non-isomorphic C$^*$-algebras with identical $\mathrm{K}$-theory and traces is 
a new phenomenon, and so our theorem constitutes a considerable strengthening of \cite[Theorem 1.1]{To2}. 

Our main theorem confirms \cite[Conjecture 6.10]{To3}.  When this conjecture was made, there were
two main obstacles to its confirmation.  The first was the lack of a good formula for the radius of
comparison of a commutative C$^*$-algebra, a problem since overcome (see \cite{To4}).  The second
was the lack of a lower semicontinuity result (with respect to inductive limits) for the radius of comparison.  
We provide such a result in the sequel;  this is the main technical innovation of the paper.

Our paper is organised as follows:  in Section 2 we review Cuntz comparison, dimension functions, and
the radius of comparison;  Section 3 studies the behaviour of the radius of comparison with respect to inductive
limits;  Section 4 contains the proof of the main theorem.  

\vspace{3mm}
\noindent
{\it Acknowledgements.} 
We thank Wilhelm Winter both for general comments on
earlier drafts and for his help with the proof of Lemma \ref{lsclem}.

\section{Preliminaries}\label{prelim}

Let $A$ be a C$^*$-algebra, and let $\mathrm{M}_n(A)$ denote the $n \times n$ 
matrices whose entries are elements of $A$.  If $A = \mathbb{C}$, then we simply write $\mathrm{M}_n$.
Let $\mathrm{M}_{\infty}(A)$ denote the algebraic limit of the
direct system $(\mathrm{M}_n(A),\phi_n)$, where $\phi_n:\mathrm{M}_n(A) \to \mathrm{M}_{n+1}(A)$
is given by
\[
a \mapsto \left( \begin{array}{cc} a & 0 \\ 0 & 0 \end{array} \right).
\]
Let $\mathrm{M}_{\infty}(A)_+$ (resp. $\mathrm{M}_n(A)_+$)
denote the positive elements in $\mathrm{M}_{\infty}(A)$ (resp. $\mathrm{M}_n(A)$). 
Given $a,b \in \mathrm{M}_{\infty}(A)_+$, we say that $a$ is {\it Cuntz subequivalent} to
$b$ (written $a \precsim b$) if there is a sequence $(v_n)_{n=1}^{\infty}$ of
elements of $\mathrm{M}_{\infty}(A)$ such that
\[
||v_nbv_n^*-a|| \stackrel{n \to \infty}{\longrightarrow} 0.
\]
We say that $a$ and $b$ are {\it Cuntz equivalent} (written $a \sim b$) if
$a \precsim b$ and $b \precsim a$.  This relation is an equivalence relation,
and we write $\langle a \rangle$ for the equivalence class of $a$.  The set
\[
W(A) := \mathrm{M}_{\infty}(A)_+/ \sim
\] 
becomes a positively ordered Abelian monoid when equipped with the operation
\[
\langle a \rangle + \langle b \rangle = \langle a \oplus b \rangle
\]
and the partial order
\[
\langle a \rangle \leq \langle b \rangle \Leftrightarrow a \precsim b.
\]
In the sequel, we refer to this object as the {\it Cuntz semigroup} of $A$.  The
Grothendieck enveloping group of $W(A)$ is denoted $\mathrm{K}_0^*(A)$.

Given $a \in \mathrm{M}_{\infty}(A)_+$ and $\epsilon > 0$, we denote by 
$(a-\epsilon)_+$ the element of $C^*(a)$ corresponding (via the functional
calculus) to the function
\[
f(t) = \mathrm{max}\{0,t-\epsilon\}, \ t \in \sigma(a).
\]
(Here $\sigma(a)$ denotes the spectrum of $a$.)  
The proposition below collects some facts about Cuntz subequivalence due 
to Kirchberg and R{\o}rdam.

\begin{props}[Kirchberg-R{\o}rdam (\cite{KR}), R{\o}rdam (\cite{R4})]\label{basics}
Let $A$ be a C$^*$-algebra, and $a,b \in A_+$.  
\begin{enumerate}
\item[(i)] $(a-\epsilon)_+ \precsim a$ for every $\epsilon > 0$.
\item[(ii)] The following are equivalent:
\begin{enumerate}
\item[(a)] $a \precsim b$;
\item[(b)] for all $\epsilon > 0$, $(a-\epsilon)_+ \precsim b$;
\item[(c)] for all $\epsilon > 0$, there exists $\delta > 0$ such that $(a-\epsilon)_+ \precsim (b-\delta)_+$.
\end{enumerate}
\item[(iii)] If $\epsilon>0$ and $||a-b||<\epsilon$, then $(a-\epsilon)_+ \precsim b$.
\end{enumerate}
\end{props} 

Now suppose that $A$ is unital and stably finite, and denote by $\mathrm{QT}(A)$
the space of normalised 2-quasitraces on $A$ (v. \cite[Definition II.1.1]{BH}).
Let $S(W(A))$ denote the set of additive and order preserving maps $s$ from $W(A)$ to $\mathbb{R}^+$
having the property that $s(\langle 1_A \rangle) = 1$.
Such maps are called {\it states}.  Given $\tau \in \mathrm{QT}(A)$, one may 
define a map $s_{\tau}:\mathrm{M}_{\infty}(A)_+ \to \mathbb{R}^+$ by
\begin{equation}\label{ldf}
s_{\tau}(a) = \lim_{n \to \infty} \tau(a^{1/n}).
\end{equation}
This map is lower semicontinous, and depends only on the Cuntz equivalence class
of $a$.  It moreover has the following properties:
\vspace{2mm}
\begin{enumerate}
\item[(i)] if $a \precsim b$, then $s_{\tau}(a) \leq s_{\tau}(b)$;
\item[(ii)] if $a$ and $b$ are mutually orthogonal, then $s_{\tau}(a+b) = s_{\tau}(a)+s_{\tau}(b)$;
\item[(iii)] $s_{\tau}((a-\epsilon)_+) \nearrow s_{\tau}(a)$ as $\epsilon \to 0$. 
\end{enumerate}
\vspace{2mm}
Thus, $s_{\tau}$ defines a state on $W(A)$.
Such states are called {\it lower semicontinuous dimension functions}, and the set of them 
is denoted $\mathrm{LDF}(A)$. $\mathrm{QT}(A)$ is a simplex (\cite[Theorem II.4.4]{BH}), 
and the map from $\mathrm{QT}(A)$ to $\mathrm{LDF}(A)$ defined by (\ref{ldf}) is 
bijective and affine (\cite[Theorem II.2.2]{BH}).  A {\it dimension function} on
$A$ is a state on $\mathrm{K}_0^*(A)$, assuming that the latter has been equipped
with the usual order coming from the Grothendieck map.  The set of dimension functions
is denoted $\mathrm{DF}(A)$.  $\mathrm{LDF}(A)$ is a (generally proper) face of
$\mathrm{DF}(A)$.  If $A$ has the property that $a \precsim b$ whenever $s(a) < s(b)$
for every $s \in \mathrm{LDF}(A)$, then we say that $A$ has 
{\it strict comparison of positive elements}.

Finally, we recall the definition of the radius of comparison.
\begin{dfs}[Definition 6.1, \cite{To3}]\label{rc}
Say that $A$ has $r$-comparison if whenever one has positive elements
$a,b \in \mathrm{M}_{\infty}(A)$ such that
\[
s(\langle a \rangle) + r < s(\langle b \rangle), \ \forall s \in \mathrm{LDF}(A),
\]
then $\langle a \rangle \leq \langle b \rangle$ in $W(A)$.
Define the radius of comparison
of $A$, denoted $\mathrm{rc}(A)$, to be 
\[
\mathrm{inf} \{r \in \mathbb{R}^+ | \ (W(A), \langle 1_A \rangle) \ \mathrm{has} \ r-\mathrm{comparison} \ \}
\]
if it exists, and $\infty$ otherwise.
\end{dfs}

\section{The radius of comparison and inductive limits}\label{rc}

\noindent

In this section we prove that the radius of comparison is
lower semicontinuous with respect to inductive limits provided that 
the limit algebra is simple and the inductive sequence is
an ``efficient decomposition'' with respect to traces.  We then show that if, as conjectured by
Blackadar and Handelman, one has $\mathrm{LDF}(A)$ weak-$*$ dense in $\mathrm{DF}(A)$
for unital and stably finite $A$, then this ``efficient decomposition'' assumption
is unneccessary.

Our first lemma is a state-theoretic analogue of the implication (i) $\Rightarrow$ (iii) of Proposition \ref{basics}.

\begin{lms}\label{stateapprox}
Let $A$ be a simple, unital, and stably finite C$^*$-algebra.  Let $a,b \in \mathrm{M}_{\infty}(A)_+$
and $r > 0$ be such that 
\[
s(a) + r < s(b), \ \forall s \in \mathrm{LDF}(A).
\]
Then, given $\epsilon > 0$, there exists $\delta>0$ such that
\[
s((a-\epsilon)_+) + r < s((b-\delta)_+), \ \forall s \in \mathrm{LDF}(A).
\]
\end{lms}

\begin{proof}
Write $s_{\tau}$ for the element of $\mathrm{LDF}(A)$ induced by $\tau \in \mathrm{QT}(A)$,
and suppose first that $\langle a \rangle = \langle p \rangle$ for some projection
$p \in \mathrm{M}_{\infty}(A)$.  Then, the map $\tau \mapsto s_{\tau}(a)$ defines a continuous
function on $\mathrm{QT}(A)$.  By \cite[Proposition 2.7]{PT}, the map $\tau \mapsto s_{\tau}(b)$
is lower semicontinuous on $\mathrm{QT}(A)$ for any $b \in \mathrm{M}_{\infty}(A)_+$.
Thus, the map $\tau \mapsto s_{\tau}(b) - s_{\tau}(a) - r$ is lower semicontinuous and strictly
positive, and so achieves a minimum value $\eta > 0$ on the compact set $\mathrm{QT}(A)$.
For each $\tau \in \mathrm{QT}(A)$ there is an open neighbourhood $U_{\tau}$ of $\tau$
such that
\[
|s_{\gamma}(a) - s_{\tau}(a)| < \eta/3, \ \forall \gamma \in U_{\tau}.
\]

We have that $s((b-\delta)_+) \nearrow s(b)$ as $\delta \to 0$, for every $s \in \mathrm{LDF}(A)$.
Choose, then, for each $\tau \in \mathrm{QT}(A)$, a $\delta_{\tau}>0$ such that 
\[
s_{\tau}(b) - s_{\tau}((b-\delta_{\tau})_+) < \eta/3.
\]
The lower semicontinuity of $\gamma \mapsto s_{\gamma}((b-\delta_{\tau})_+)$ 
ensures that for each $\tau \in \mathrm{QT}(A)$,
there is an open neighbourhood $V_{\tau}$ of $\tau$ such that
\[
s_{\gamma}((b-\delta_{\tau})_+) > s_{\tau}((b-\delta_{\tau})_+) - \eta/3, \ \forall \gamma \in V_{\tau}.
\]
Put $W_{\tau} = U_{\tau} \cap V_{\tau}$.  Then,
\[
s_{\gamma}((b-\delta_{\tau})_+) > s_{\gamma}(a) + r \geq s_{\gamma}((a-\epsilon)_+) + r, \ \forall \gamma \in W_{\tau}.
\]
The $W_{\tau}$ cover the compact set $\mathrm{QT}(A)$ and so there exist $\tau_1,\ldots,\tau_n \in \mathrm{QT}(A)$
such that
\[
\mathrm{QT}(A) = W_{\tau_1} \cup \cdots \cup W_{\tau_n}.
\]
Put $\delta = \mathrm{min} \{\delta_{\tau_1},\ldots,\delta_{\tau_n}\}$.  Then, since 
\[
s_{\gamma}((b-\delta)_+) \geq s_{\gamma}((b-\delta_{\tau_i})_+), \ \forall \gamma \in \mathrm{QT}(A), \ 1 \leq i \leq n,
\]
we have 
\[
s((a-\epsilon)_+) + r < s((b-\delta)_+), \ \forall s \in \mathrm{LDF}(A),
\]
as desired.

Now suppose that $\langle a \rangle \neq \langle p \rangle$ for any projection
$p \in \mathrm{M}_{\infty}(A)$.  Then, by a functional calculus argument, $0$ is not
an isolated point of the spectrum $\sigma(a)$ of $a$.  The simplicity of $A$ implies
that each quasitrace on $A$ is faithful, so for each $\epsilon>0$ we have 
\[
s((a-\epsilon)_+) < s(a), \ \forall s \in \mathrm{LDF}(A).
\]
Each $\tau \in \mathrm{QT}(A)$ is implemented on $C^*(a)$ by
a probability measure $\mu_{\tau}$ on $\sigma(a)$.  Put
\[
A_{\tau} = \{\gamma \in \mathrm{QT}(A) \ | \ \mu_{\gamma}([\epsilon,\infty) \cap \sigma(a)) \geq s_{\tau}(a) \}.
\]
If $\gamma \notin A_{\tau}$, then 
\[
s_{\gamma}((a-\epsilon)_+) = \mu_{\gamma}((\epsilon,\infty) \cap \sigma(a)) \leq \mu_\gamma([\epsilon,\infty)
\cap \sigma(a)) < s_{\tau}(a).
\]
Clearly, $\tau \notin A_{\tau}$.  Let $(\gamma_n)$ be a convergent sequence in $A_{\tau}$ with
limit $\gamma$.  This implies that $\mu_{\gamma_n} \to \mu_{\gamma}$ in measure on $\sigma(a)$.
By Portmanteau's Theorem, this implies that $\mu_{\gamma}(C) \geq \limsup \mu_{\gamma_n}(C)$
for every closed subset of $\sigma(a)$.  Since $[\epsilon,\infty) \cap \sigma(a)$ is closed we
have
\[
\mu_{\gamma}([\epsilon,\infty) \cap \sigma(a)) \geq \limsup \mu_{\gamma_n}([\epsilon,\infty) \cap \sigma(a)) \geq s_{\tau}(a),
\]
whence $\gamma \in A_{\tau}$ and $A_{\tau}$ is closed.  Thus $U_{\tau} := A_{\tau}^c$ is an open neighbourhood of
$\tau$, and 
\[
s_{\gamma}((a-\epsilon)_+) < s_{\tau}(a), \ \forall \gamma \in U_{\tau}.
\]

Find (as above) for each $\tau \in \mathrm{QT}(A)$ a $\delta_{\tau}>0$ such that 
$s_{\tau}((b-\delta_{\tau})_+) > s_{\tau}(a) +r $.  Apply the lower semicontinuity of 
$\gamma \mapsto s_{\gamma}((b-\delta_{\tau})_+)$ to find an open neighbourhood
$V_{\tau}$ of $\tau$ such that
\[
s_{\gamma}((b-\delta_{\tau})_+) > s_{\tau}(a) + r, \ \forall \gamma \in V_{\tau}.
\]
Put $W_{\tau} = U_{\tau} \cap V_{\tau}$.  Then,
\[
s_{\gamma}((a-\epsilon)_+) + r < s_{\tau}(a) +r< s_{\gamma}((b-\delta_{\tau})_+), \ \forall \gamma \in W_{\tau}.
\]
Use the compactness argument above to find $\delta>0$ such that 
\[
s((a-\epsilon)_+) + r < s((b-\delta)_+), \ \forall s \in \mathrm{LDF}(A).
\]

\end{proof}

The next lemma is due to M. R{\o}rdam.  Its proof can be found in \cite{To4}.

\begin{lms}[R{\o}rdam, \cite{R7}]\label{posapprox1}
Let $A$ be a C$^*$-algebra and $\{A_i\}_{i \in I}$ a collection of C$^*$-subalgebras whose
union is dense.   Then, for every $a \in \mathrm{M}_{\infty}(A)_+$ and $\epsilon > 0$
there exists $i \in \mathbb{N}$ and $\tilde{a} \in \mathrm{M}_{\infty}(A_i)$ such that
\[
(a-\epsilon)_+ \precsim \tilde{a} \precsim (a-\epsilon/2)_+ \precsim a 
\]
in $\mathrm{M}_{\infty}(A)$
\end{lms}

\begin{props}\label{rclsc}
Let $A = \lim_{i \to \infty}(A_i,\phi_i)$ be a simple, unital, and stably finite C$^*$-algebra with $\phi_i$ injective. 
Suppose this decomposition satisfies the following property:  given $\epsilon>0$ and
$a,b \in \mathrm{M}_{\infty}(\cup_{i=1}^{\infty} A_i)_+$ such that
\[
s(a) < s(b), \ \forall s \in \mathrm{LDF}(A),
\]
there is a $j \in \mathbb{N}$ such that
\[
s(a) < s(b)+\epsilon, \ \forall s \in \mathrm{LDF}(A_j).
\]
Then, 
\[
\mathrm{rc}(A) \leq \liminf_{i \to \infty} \mathrm{rc}(A_i).
\]
\end{props}

\begin{proof}
The theorem is trivial if 
\[
\liminf \mathrm{rc}(A_i) = \infty, 
\]
so suppose
that 
\[
r := \liminf \mathrm{rc}(A_i) < \infty.
\]
Passing to a subsequence if necessary, we assume that $(\mathrm{rc}(A_i))_{i=1}^{\infty}$
is decreasing.

Let there be given $a,b \in \mathrm{M}_n(A)_+ \hookrightarrow \mathrm{M}_{\infty}(A)_+$ 
 and $m>r$ such that
\[
s(a) + m < s(b), \ \forall s \in \mathrm{LDF}(A).
\]
By \cite[Proposition 6.3]{To3}, it will suffice to prove that $a \precsim b$.
Let $\epsilon > 0$ be given, and use Lemma \ref{stateapprox} to find 
a $\delta>0$ such that
\[
s((a-\epsilon/2)_+) + m  < s((b-\delta)_+), \ \forall s \in \mathrm{LDF}(A).
\]
Find, using Lemma \ref{posapprox1}, a positive element $\tilde{a} \in A_j$, some $j \in \mathbb{N}$, such that
\[
||(a - \epsilon/2)_+ - \tilde{a}|| < \epsilon/2
\]
and 
\[
(a-\epsilon)_+ \precsim \tilde{a} \precsim (a-\epsilon/2)_+.
\]
Put $\epsilon^{'} = \mathrm{min} \{\epsilon, \delta \}$, and find a positive element $\tilde{b}$ in some $A_j$ 
(we may assume that it is the same $A_j$ that contains $\tilde{a}$) such that $||b - \tilde{b}||< \epsilon^{'}$
and $(b-\delta)_+ \precsim \tilde{b} \precsim b$.  Finally, assume that $j \in \mathbb{N}$ has been
chosen large enough to ensure that $\mathrm{rc}(A_j) < m$.
We now have 
\[
s_{\gamma}(\tilde{a}) + m \leq s_{\gamma}((a-\epsilon/2)_+) +m < s_{\gamma}((b-\delta)_+) \leq s_{\gamma}(\tilde{b}), \ \forall \gamma \in \mathrm{QT}(A).
\]
Choose $\eta > 0$ such that $m-\eta > \mathrm{rc}(A_j)$.  By our hypothesis we may assume,
upon incresing $j$ if necessary, that
\[
s_{\gamma}(\tilde{a}) + m-\eta < s_{\gamma}(\tilde{b}), \ \forall \gamma \in \mathrm{QT}(A_j).
\]
Since $\mathrm{rc}(A_j) < m-\eta$, we conclude that
\[
(a-\epsilon)_+ \precsim \tilde{a} \precsim \tilde{b} \precsim b;
\]
$\epsilon$ was arbitrary, and the proposition follows.
\end{proof}

The reader may wonder whether the ``extra'' hypothesis of Proposition \ref{rclsc}
--- the ``efficient decomposition'' hypothesis alluded to at the beginning of this
subsection --- can be removed.  Indeed, if $a$ and $b$ as in the proposition are projections, then there is 
always a natural number $j$ such that $s_{\tau}(a) = \tau(a) < \tau(b) = s_{\tau}(b)$
for each $\tau \in \mathrm{QT}(A_j)$.  (To the best of our knowledge, this was first observed by Blackadar in \cite{Bl3}.)
We point out why this argument does not carry over to the setting of positive elements 
and lower semicontinuous dimension functions after proving the next lemma.  We
thank Wilhelm Winter for pointing out the ultrafilter argument used in the proof.

\begin{lms}\label{lsclem}
Let $A = \lim_{i \to \infty}(A_i,\phi_i)$ be simple and unital, with each $A_i$ stably finite
and each $\phi_i$ injective.  Suppose that $a,b \in \mathrm{M}_{\infty}(\cup_{i=1}^{\infty}A_i)_+$ are such that
\[
s(a) < s(b), \ \forall s \in \mathrm{DF}(A).
\]
Then, there is some $j \in \mathbb{N}$ such that $a,b \in \mathrm{M}_{\infty}(A_j)_+$
and
\[
s(a) < s(b), \ \forall s \in \mathrm{DF}(A_j).
\]
\end{lms}

\begin{proof}
Let $B$ be the C$^*$-algebra consisting of all bounded sequences
\[
(a_1,a_2,a_3,\ldots),
\]
where $a_i \in A_i$, and let $I \subseteq B$ be the closed two-sided
ideal of sequences such that $a_i \to 0$ as $i \to \infty$.  Then, there
is a $*$-monomorphism $\iota:A \to B/I$ given by
\[
a \mapsto (\underbrace{0,\ldots,0}_{i-1 \ \mathrm{times}},a,\phi_i(a),\phi_{i+1}(\phi_i(a)),\ldots) + I
\]
on $A_i \subseteq A$, $i \in \mathbb{N}$, and extended by continuity.  

We may assume, by truncating our inductive sequence if necessary, that $a,b \in A_1$.
Suppose, contrary to our desired conclusion, that for each $i \in \mathbb{N}$ there exists
a dimension function $d_i \in \mathrm{DF}(A_i)$ satisfying
\[
d_i(a) \geq d_i(b).
\]
Let $\omega$ be a free ultrafilter on $\mathbb{N}$, and let
$s$ be the map given by
\[
s(a_1,a_2,\ldots) = \lim_{\omega} d_i(a_i).
\]
It is straightforward to check that $s \in \mathrm{DF}(A)$.  But then
$s(b) \geq s(a)$, contrary to our assumption.     
\end{proof}

\noindent
The free ultrafilter approach in Lemma \ref{lsclem} can be applied after replacing 
the $d_i$ with lower semicontinuous dimension functions $s_{\tau_i} \in \mathrm{LDF}(A_i)$,
but it is then unclear whether the resulting dimension function is lower semicontinuous.
Alternatively, one can use the $\tau_i$ themselves to define, via the free ultrafilter,
a faithful trace $\tau$ on $A$, but it is then unclear whether $s_{\tau}(a) \geq s_{\tau}(b)$.
At issue is an interchanging of the limit over $\omega$ and the limit appearing in
the definition of a lower semicontinuous dimension function.  In any case, there is no
obvious way to prove the lemma upon substituting lower semicontinuous dimension functions
for dimension functions.  But this difficulty vanishes if $\mathrm{LDF}(A)$ is dense in $\mathrm{DF}(A)$
whenever $A$ is unital and stably finite.

\begin{thms}
Let $A = \lim_{i \to \infty}(A_i,\phi)$ be unital and simple, with each $A_i$ stably finite
and each $\phi_i$ injective.  Also suppose that $\mathrm{LDF}(A)$ is dense in $\mathrm{DF}(A)$.
Then,
\[
\mathrm{rc}(A) \leq \liminf_{i \to \infty} \mathrm{rc}(A_i).
\]
\end{thms}

\begin{proof}  As in the proof of Proposition \ref{rclsc}, we assume that
\[
r:=\liminf_{i \to \infty} \mathrm{rc}(A_i) < \infty.
\]
Let $\epsilon>0$ and $a,b \in A_+$ be given.  For any $d \in \mathrm{DF}(A)$, there is a sequence
$(s_{\tau_i})_{i=1}^{\infty}$ in $\mathrm{LDF}(A)$ converging to $d$.  It 
follows that $d(a) \leq d(b) < d(b) +\epsilon$ whenever 
\[
s_{\tau}(a) < s_{\tau}(b), \forall s_{\tau} \in \mathrm{LDF}(A).
\]
By Lemma \ref{lsclem} there is some $i \in \mathbb{N}$ such that 
$d(a) < d(b) + \epsilon$ for every $d \in \mathrm{DF}(A_i) \supseteq 
\mathrm{LDF}(A_i)$.  Thus, $A$ satisfies the hypotheses of Proposition
\ref{rclsc}, and the theorem follows.
\end{proof}

\section{The main result}\label{ncp}

\begin{thms}\label{ce}
There is a simple, separable, and nuclear C$^*$-algebra $A$ such that
for any natural numbers $n \neq m$ one has:
\begin{enumerate}
\item[(i)] $\mathrm{M}_n(A) \ncong \mathrm{M}_{m}(A)$;
\item[(ii)] $\mathrm{M}_n(A)$ and $\mathrm{M}_{m}(A)$ agree on the Elliott invariant;
\item[(iii)] $\mathrm{M}_n(A)$ and $\mathrm{M}_{m}(A)$ are shape equivalent;
\item[(iv)] $\mathrm{M}_n(A)$ is non-$\mathcal{Z}$-stable, and has stable 
rank one, real rank one, and property (SP);
\item[(v)] $(V(\mathrm{M}_n(A)),[1_{\mathrm{M}_n(A)}]) \cong (\mathbb{Q}^+,1)$;
in particular, $\mathrm{K}_0(\mathrm{M}_n(A))$ is a divisible and weakly unperforated partially ordered group.
\end{enumerate}
\end{thms}

Before proceeding with the proof we briefly recall some terminology.
Let $X$ and $Y$ be compact Hausdorff spaces, and let $m,n \in \mathbb{N}$ be such
that $m|n$.  Recall that a $*$-homomorphism
\[
\phi:\mathrm{M}_m(\mathrm{C}(X)) \to \mathrm{M}_{n}(\mathrm{C}(Y))
\]
is called {\it diagonal} if 
\[
\phi(f) = \bigoplus_{i=1}^{n/m} f \circ \lambda_i,
\]
where each $\lambda_i:Y \to X$ is continuous.  The $\lambda_i$ are called {\it eigenvalue maps}.

\begin{proof}
$A$ will be constructed as per the general framework set out in \cite{V1},
and will be identical to the construction in the proof of \cite[Theorem 1.1]{To2}.
It is necessary, however, to recall the details of the construction, as they are
essential to proving that $\mathrm{rc}(A)$ is finite and non-zero.
 
Put $X = [-1,1]^3$.   Put $X_1 = X \times X$, and put $X_{i+1} = (X_i)^{n_i}$ --- the $n_i$-fold
Cartesian product of $X_i$ with itself.  Let $\pi_i^j:X_{i+1} \to X_i$, $1 \leq j \leq n_i$ be the co-ordinate
projections.  Let $A_i$ be the homogeneous C$^*$-algebra $\mathrm{M}_{m_i} \otimes 
\mathrm{C}(X_i)$, where $m_i$ is a natural number to be specified, and let $\phi_i:A_i \to A_{i+1}$
be the $*$-homomorphism given by 
\[
\phi_i(a)(x) = \mathrm{diag}\left( a \circ \pi_i^1(x),\ldots,a \circ \pi_i^{n_i}(x),a(x_i^1),\ldots,a(x_i^i) \right), \ \forall x \in X_{i+1},
\]
where $x_i^1,\ldots,x_i^i \in X_i$ are to be specified.  Put $A = \lim_{i \to \infty} (A_i,\phi_i)$,
and define
\[
\phi_{i,j} := \phi_{j-1} \circ \cdots \circ \phi_i.
\]
Let $\phi_{i\infty}:A_i \to A$ be the canonical map.
Assume that the $n_i$ have been chosen so that
\[
\prod_{i=1}^{\infty} \frac{n_i}{i+n_i} \neq 0,
\] 
and that the $x_i^1,\ldots,x_i^i$ have been chosen to ensure that $A$ is simple (this can always be arranged
--- cf. \cite{V1}).  Finally, assume that our choice of the $n_i$ is such that every natural number
divides some $m_i$.  Put $m_1=2$.   

\vspace{3mm}
\noindent
{\bf (i).} We will first prove that the inductive sequence $(A_i,\phi_i)$ satisfies the hypotheses of Proposition \ref{rclsc},
i.e., that given $\epsilon>0$ and
$a,b \in \mathrm{M}_{\infty}(\cup_{i=1}^{\infty} A_i)_+$ such that
\[
s(a) < s(b), \ \forall s \in \mathrm{LDF}(A),
\]
there is an $i \in \mathbb{N}$ such that
\[
s(a) < s(b)+\epsilon, \ \forall s \in \mathrm{LDF}(A_i).
\]
We will accomplish this by finding for each $\epsilon>0$ an $i \in \mathbb{N}$ such that
the following holds:  for any $\tau \in \mathrm{T}(A_i)$ there is some $\gamma \in \mathrm{T}(A)$ satisfying
\[
\phi_{i\infty}^{\sharp}(\gamma) = (1-\lambda) \tau + \lambda \eta
\]
for some $\lambda$ such that $0 < \frac{\lambda}{1-\lambda} < \epsilon$ and $\eta \in \mathrm{T}(A_i)$.  To see that
this will suffice, take any $\tau \in \mathrm{T}(A_i)$
and $a,b \in \mathrm{M}_{\infty}(\cup_{i=1}^{\infty} A_i)_+$ such that
\[
s(a) < s(b), \ \forall s \in \mathrm{LDF}(A).
\]
Then, we have
\begin{eqnarray*}
s_{\tau}(a) & \leq & s_{\tau}(a) + \frac{\lambda}{1-\lambda} s_{\eta}(a) \\
& = & \frac{1}{1-\lambda} s_{(1-\lambda)\tau+ \lambda \eta}(a) \\
& = & \frac{1}{1-\lambda} s_{\phi_{i\infty}^{\sharp}(\gamma)}(a) \\
& = & \frac{1}{1-\lambda} s_{\gamma}(a) \\
& < & \frac{1}{1-\lambda} s_{\gamma}(b) \\
& = & \frac{1}{1-\lambda} s_{\phi_{i\infty}^{\sharp}(\gamma)}(b) \\
& = & s_{\tau}(b) + \frac{\lambda}{1-\lambda} s_{\eta}(b) \\
& < & s_{\tau}(b) + \epsilon.
\end{eqnarray*}

Let $\epsilon > 0$ be given.  For any $i \in \mathbb{N}$, $\tau \in \mathrm{T}(A_i)$, and $j >i$,
let $\tau_{i,j} \in \mathrm{T}(A_j)$ be the trace corresponding to the product measure
\[
\underbrace{\tau \times \cdots \times \tau}_{n_{i+1} n_{i+2} \cdots n_j \ \mathrm{times}}.
\]
If $i \leq k < j$, then define $\tau_k^j:= \phi_{k,j}^{\sharp}(\tau_{i,j})$.  
A staightforward calculation shows that 
$\tau_k^j = (1-\lambda_k^j) \tau_{i,k} + \lambda_k^j \eta_k^j$ for some $\eta_k^j \in \mathrm{T}(A_k)$
and $0 < \lambda_k^j < 1$.  In fact, if $N_{k,j}$ denotes the number of eigenvalue maps of $\phi_{k,j}$
which are co-ordinate projections onto $X_k$, then 
\[
\lambda_k^j = \frac{N_{k,j}}{\mathrm{mult}(\phi_{k,j})}.
\]
By construction $N_{k,j}/\mathrm{mult}(\phi_{k,j})$ is a decreasing sequence with strictly positive
limit $1-\lambda_k$.  By increasing $i$ if necessary, we may assume that $\lambda_i/(1-\lambda_i) < \epsilon$.
We claim that for each $k \geq i$, 
\[
(1-\lambda_k^j) \tau_{i,k} + \lambda_k^j \eta_k^j \stackrel{j \to \infty}{\longrightarrow} 
(1-\lambda_k) \tau_{i,k} + \lambda_{k} \eta_k
\]
for some $\eta_k \in \mathrm{T}(A_k)$.  It will suffice to prove that $\eta_k^j$
is a Cauchy sequence in $j$.
Let $R_{k,j}$ be the multiset whose elements are the points of $X_k$ which appear as point
evaluations in $\phi_{k,j}$.  We have that
\[
\eta_k^j = \frac{1}{|R_{k,j}|} \sum_{x \in R_{k,j}} ev_x;
\]
$\eta_k^j$ is a finite convex combination of extreme traces on $A_k$.  $R_{k, j+1}$ is formed
by taking the union (with multiplicity) of $\mathrm{mult}(\phi_{j,j+1})$ copies of $R_{k,j}$
and some other multiset $S_{j+1}$.  Thus,
\begin{eqnarray*}
\eta_k^{j+1} & = & \frac{1}{|R_{k,j+1}|} \sum_{x \in R_{k,j+1}} ev_x \\ 
& = & \frac{\mathrm{mult}(\phi_{j,j+1})}{\mathrm{mult}(\phi_{j,j+1})|R_{j,k}| + |S_{j+1}|} \sum_{x \in R_{k,j}} ev_x \\
& & + \frac{1}{\mathrm{mult}(\phi_{j,j+1})|R_{j,k}| + |S_{j+1}|} \sum_{x \in S_{j+1}} ev_x \\
& = & \frac{\mathrm{mult}(\phi_{j,j+1})|R_{j,k}|}{\mathrm{mult}(\phi_{j,j+1})|R_{j,k}| + |S_{j+1}|} \eta_k^j \\
& & + \frac{1}{\mathrm{mult}(\phi_{j,j+1})|R_{jk}| + |S_{j+1}|} \sum_{x \in S_{j+1}} ev_x.
\end{eqnarray*}
It follows that for any normalised element $f \in A_k$, we have that
\[
|\eta_k^j(f) - \eta_k^{j+1}(f)| \leq \frac{2|S_{j+1}|}{\mathrm{mult}(\phi_{j,j+1})|R_{j,k}| + |S_{j+1}|}.
\]
By construction, the right hand side vanishes as $j \to \infty$, proving the claim.

We now have that for each $k \in \mathbb{N}$, the sequence $\tau_k^j$ converges as $j \to \infty$.
Call the limit $\tau_k$.  Since $\tau_k^j = \phi_{k,k+1}^{\sharp}(\tau_{k+1}^j)$, we have
that $\tau_k = \phi_{k,k+1}^{\sharp}(\tau_{k+1})$ for every $k \geq i$.  It follows that
the sequence $(\tau_i,\tau_{i+1},\ldots)$ defines a point in limit of the inverse
system $(\mathrm{T}(A_i),\phi_i^{\sharp})$, i.e., a point in $\mathrm{T}(A)$.
Thus, 
\[
\tau_i = (1-\lambda_i) \tau + \lambda_{i} \eta_i
\]
is the image of some $\gamma \in \mathrm{T}(A)$ under the map $\phi_{i\infty}^{\sharp}$.
Since $0 < \lambda_i/(1-\lambda_i) < \epsilon$, we have established the hypotheses of 
Proposition \ref{rclsc} for the inductive system $(A_i,\phi_i)$.  
It is clear from our construction and \cite[Theorem 4.2]{To4} that $\mathrm{rc}(A_i) < 10$, $\forall i \in \mathbb{N}$.
Proposition \ref{rclsc} then shows that $\mathrm{rc}(A) < \infty$. 
The proof of \cite[Theorem 1.1]{To2} shows that $\mathrm{rc}(A) > 0$, so that $\mathrm{rc}(A)$ is finite and nonzero.

Now let $m \neq n$ be natural numbers.  \cite[Proposition 6.2, (i)]{To3} shows that 
\[
\mathrm{rc}(\mathrm{M}_n(A)) = \mathrm{rc}(A)/n \neq \mathrm{rc}(A)/m = \mathrm{rc}(\mathrm{M}_m(A)),
\]
whence $\mathrm{M}_n(A) \ncong \mathrm{M}_m(A)$, as desired.

\vspace{3mm}
\noindent
{\bf (ii).}
 By the contractibility of $X_i$ we have $\mathrm{K}_0(A_i) \cong \mathbb{Z}$, $\mathrm{K}_1(A_i) = 0$,
and $\mathrm{K}_0(\phi_i)(1) = m_i$.  It follows from our assumption on the $m_i$s that 
\[
(\mathrm{K}_0(A),\mathrm{K}_0(A)^+,[1_A]) \cong (\mathbb{Q},\mathbb{Q}^+,1),
\]
and this same isomorphism clearly holds for any matrix algebra over $A$. 

Since $\mathrm{K}_0(\mathrm{M}_n(A)) \cong \mathbb{Q}$, there is a unique pairing
between traces and $\mathrm{K}_0$.  Thus, all of the non-stable information in the Elliott
invariant of $\mathrm{M}_n(A)$ is independent of $n$.  The remaining
elements of the Elliott invariant are stable isomorphism invariants, and are thus
also independent of $n$.   

\vspace{3mm}
\noindent
{\bf (iii).}
We will prove that $A$ and $\mathrm{M}_2(A)$ are shape equivalent.  The argument for $\mathrm{M}_n(A)$ 
and $\mathrm{M}_m(A)$ is similar.

By compressing the inductive sequence $(A_i,\phi_i)$ if necessary, we may assume that $2m_i|m_{i+1}$
for all $i \in \mathbb{N}$.
To prove that $A$ and $\mathrm{M}_2(A)$ are shape equivalent, it will suffice to
find sequences of $*$-homomorphisms $(\eta_i)_{i=1}^{\infty}$ and $(\mu_i)_{i=1}^{\infty}$ making the diagram
\[
\xymatrix{
{A_1}\ar[rr]^-{\phi_1}\ar[dd]^{\eta_1}&&{A_2}\ar[rr]^-{\phi_2}\ar[dd]^{\eta_2}&&{A_3}\ar[r]\ar[dd]^{\eta_3}& {\cdots}\ar[r] & A \\
&&&&&& \\
{\mathrm{M}_2 \otimes A_1}\ar[rr]^-{\mathrm{id} \otimes \phi_1}\ar[uurr]_-{\mu_1}&&
{\mathrm{M}_2 \otimes A_2}\ar[rr]^-{\mathrm{id} \otimes \phi_2}\ar[uurr]_-{\mu_2}&&
{\mathrm{M}_2 \otimes A_3}\ar[r]\ar[ur]_-{\mu_3}& {\cdots}\ar[r] & \mathrm{M}_2 \otimes A
}
\]
commute up to homotopy.  Fix a point $y_i \in X_i$ for every $i \in \mathbb{N}$, and
define:
\[
\eta_i(f) = f(y_i) \oplus f(y_i); \ \  \mu_i(g) = \bigoplus_{l=1}^{\frac{n_{i+1}}{2 n_i}} g(y_i).
\]
Then, for every $i \in \mathbb{N}$, the maps $\phi_i, \mathrm{id}_{\mathrm{M}_2} \otimes \phi_i,
\mu_i \circ \eta_i$, and $\eta_i \circ \mu_{i-1}$ are diagonal; they have the form
\[
f \mapsto \bigoplus_{j=1}^{k} f \circ \lambda_j,
\]
where the $\lambda_j$ are continuous maps from the spectrum of the target algebra to the spectrum
of the source algebra.  The latter are always pairwise homotopic in our setting, since each $X_i$
is contractible.  It follows that any two diagonal maps from $A_i$ to $A_{i+1}$, or from 
$\mathrm{M}_2(A_i)$ to $\mathrm{M}_2(A_{i+1})$, are homotopic, and so our diagram commutes
up to homotopy as required.  

\vspace{3mm}
\noindent
{\bf (iv).}
That the stable rank and the real rank of $A$ are one follows from \cite[Proposition 10]{V1} and the discussion
following it.  Both the stable and real rank are monotone decreasing under taking tensor products
with full matrix algebras, with the caveat that in both cases an algebra with non-minimal rank is not Morita
equivalent to an algebra with minimal rank (\cite[Theorem 3.6]{Ri1} and \cite{BP}).
It follows that every matrix algebra over $A$ has stable rank one and real rank one.  

Villadsen
proves in \cite{V1} that his choice of points $x_i^1,\ldots,x_i^i \in X_i$, $i \in \mathbb{N}$, 
ensures that for every open $V \subseteq X_j$, there is a $k>j$ such that the composed map
\[
\phi_{jk} := \phi_{k-1} \circ \cdots \circ \phi_j
\]
has an eigenvalue map which is a point evaluation with range in $V$.  It follows that the 
hereditary subalgebra generated by the image $\phi_{i \infty}(a)$ of some nonzero positive $a \in A_i$
contains a nonzero projection.  Now let $a \in A$ be positive and nonzero.  Choose $\epsilon>0$ such
that $\epsilon \ll ||a||$, and 
find a positive element $\tilde{a}$ in some $A_i$ such that $||a - \tilde{a}||< \epsilon$.
Then, $(\tilde{a}-\epsilon)_+ \precsim a$, and $(\tilde{a}-\epsilon)_+$ is a nonzero positive
element of $A_i$.  Let $p$ be a nonzero projection contained in the hereditary subalgebra of $A$
generated by $(\tilde{a}-\epsilon)_+$.  Then, $p \precsim a$ by transitivity.  It follows from the
proof of \cite[Proposition 2.2]{PT} that $\overline{aAa}$ contains a nonzero projection, whence
$A$ and matrix algebras over it all have property (SP).   

Every simple, unital, exact, finite, and $\mathcal{Z}$-stable C$^*$-algebra
has $\mathrm{rc}=0$ by \cite[Corollary 4.6]{R4}, whence $\mathrm{M}_n(A)$ is
not $\mathcal{Z}$-stable for any $n \in \mathbb{N}$.

\vspace{3mm}
\noindent
{\bf (v).} That 
\[
(V(\mathrm{M}_n(A)),[1_{\mathrm{M}_n(A)}]) \cong (\mathbb{Q}^+,1)
\]
follows from our calculation of $\mathrm{K}_0(\mathrm{M}_n(A))$ in (ii) and the fact that $A$ has stable
rank one. 

\end{proof}

\vspace{5mm}

\noindent
\emph{Andrew S. Toms \newline
Department of Mathematics and Statistics \newline
\hspace*{2mm} York University \newline
\hspace*{2mm} 4700 Keele St. \newline
\hspace*{2mm} Toronto, Ontario, M3J 1P3 \newline
Canada \newline}

\noindent
atoms@mathstat.yorku.ca

\end{document}